\documentclass[11pt]{elsarticle}
\usepackage[a4paper]{geometry}
\usepackage{amsfonts,amsmath,amssymb,amsthm,bm}
\usepackage{hyperref}

\newcommand{\Aa}{\mathcal{A}}
\newcommand{\At}{\mathtt{A}}
\newcommand{\Bt}{\mathtt{B}}
\newcommand{\Ct}{\mathtt{C}}
\newcommand{\Dt}{\mathtt{D}}
\newcommand{\dt}{\partial_t}
\newcommand{\eps}{\varepsilon}
\newcommand{\Ht}{\mathtt{H}}

\newcommand{\Mm}{\mathcal{M}}

\newcommand{\Ot}{\mathtt{O}}

\newcommand{\Rr}{\mathbb{R}}
\newcommand{\St}{\mathtt{S}}

\begin{document}
\begin{frontmatter}

\title{A nested Schur complement solver with
       mesh-independent convergence for the time domain
       photonics modeling}

\author[kiam]{M.A.~Botchev}
\ead{botchev@ya.ru}

\address[kiam]{Keldysh Institute of Applied Mathematics,
Russian Academy of Sciences,\\Miusskaya~Sq.~4, 125047 Moscow,
Russia}

\begin{abstract}
A nested Schur complement solver is proposed for
iterative solution of linear systems arising in 
exponential and implicit time integration of
the Maxwell equations with perfectly matched layer (PML)
nonreflecting boundary conditions.  
These linear systems are the so-called double saddle point systems
whose structure is handled by the Schur complement solver
in a nested, two-level fashion.
The solver is demonstrated to have a mesh-independent convergence
at the outer level, whereas the inner level
system is of elliptic type and thus can be treated efficiently 
by a variety of solvers.
\end{abstract}

\end{frontmatter}

\section{Introduction}
Numerical solution of the time-dependent Maxwell equations
is an important computational problem arising in various
scientific and engineering fields such as photonic crystal
modeling, gas and oil industry, biomedical simulations,
and astrophysics.
Rather often the application environment suggests that the Maxwell
equations have to be solved many times, for instance, for different
source functions or different medium parameters)~\cite{BotchevHanseUppu2018}.
The size of the spatial computational domain
as well as the necessity to solve the equations many times make this
task very demanding in terms of computational costs.
Therefore, advanced computational techniques have to be applied,
such as modern finite element 
discretizations~\cite{Descombes_ea2013_DGTD,Sarmany_ea2013} in space
and efficient integration schemes in time.
Along with multirate and implicit time integration 
schemes~\cite{VerwerBotchev09,Descombes_ea2016},
exponential time integration schemes~\cite{HochbruckOstermann2010} 
have recently been shown
promising for solving the Maxwell equations 
in time~\cite{Hochbruck_Pazur_ea2015,Botchev2016,BotchevHanseUppu2018}.

Exponential time integration schemes, which are essentially based on the
notion of the matrix exponential, are attractive not only because of their
excellent stability and accuracy properties but also due to their
efficiency and potential for a parallelism in time~\cite{PARAEXP,Kooij_ea2017}.
Most frequently, especially when the
spatially discretized Maxwell operator $\Aa$ (defined in~\eqref{mxw1} below)
is not a skew-symmetric matrix, the actions of the matrix exponential
in exponential time integration schemes are evaluated by Krylov subspace
methods.
To be efficient, Krylov subspace methods often need to rely on rational
approximations~\cite{DruskinKnizh98,MoretNovati04,EshofHochbruck06,PhD_Guettel} 
(so that the Krylov subspace is built up for
a rational function of $\Aa$ rather then for $\Aa$ itself) and
on the so-called restarting 
techniques~\cite{PhD_Niehoff,TalEzer2007,Afanasjew_ea08,PhD_Guettel,Eiermann_ea2011} 
(to keep the Krylov subspace dimension restricted).
A popular variant of the rational Krylov subspace methods is the
shift-and-invert (SAI) method~\cite{MoretNovati04,EshofHochbruck06}.
Rational Krylov subspace methods and, in particular, the SAI Krylov subspace method
as well as implicit time integration schemes involve solution
of linear systems with the matrix $I+\gamma\Aa$, with $\gamma>0$
being a given parameter (which is, in case of implicit
time stepping, the time step size).

Despite the significant progress achieved last decennia in sparse
direct linear system solvers, for three-dimensional (3D) problems
iterative linear system solvers remain the methods of choice.
The task of solving linear systems with the matrix $I+\gamma\Aa$
when $\Aa$ is a spatially discretized Maxwell operator
is especially challenging for the Maxwell equations.
This is caused not only by the fact that the matrix $\Aa$ has a saddle point
structure but also due to special nonreflecting boundary conditions
which are often imposed for the Maxwell equations.
In this paper a preconditioner is proposed to solve iteratively linear systems
with the matrix $\Aa$ when the popular perfectly matching layers (PML)
boundary conditions are imposed.  In this case the matrix $\Aa$
has the so-called double saddle point structure~\cite{BeikBenzi2018}.

This paper is organized as follows.  The problem is set up
in Section~2.  In Section~3 we present the nested Schur complement
solver.  Other possible preconditioners for these problems are discussed 
Section~4.  Section~5 is devoted to numerical experiments.
Finally, 
the conclusions are drawn in the last section whereas
some background material is given in two appendices.
 
\section{Problem formulation}
We are interested in solving a system of time-dependent three-dimensional (3D) Maxwell
equations
\begin{equation}
\label{mxw}
\left\{
\begin{alignedat}{4}
&\mu \dt  H \,=\,& - \sigma_1     &H  & \,-\, \nabla\times &E & \,+\, &J_H, \\
&\eps \dt E \,=\,&   \nabla\times &H  & \,-\, \sigma_2     &E & \,+\, &J_E, \\
\end{alignedat}
\right.  
\end{equation}
where $H=H(x,y,z,t)$ and $E=E(x,y,z,t)$ are respectively magnetic and electric fields,
$\mu$ is the magnetic permeability (as typical for photonics and gas-and-oil exploration
applications, $\mu\equiv 1$ for all the tests considered in this paper; 
however, in general one can have $\mu=\mu(x,y,z)$) 
and $\eps=\eps(x,y,z)>0$ is the electric permittivity.
Furthermore, $\sigma_{1,2}=\sigma_{1,2}(x,y,z)\geqslant 0$ are the conduction terms, 
such that $\sigma_2$
contains real physical conductivity as well as additional artificial conductivity
related to nonreflective boundary conditions (in this work we use
stretched coordinate formulation of the 
perfectly matched layers, PML, boundary conditions~\cite{Johnson2010_PML}), 
whereas $\sigma_1$ normally contains artificial PML conductivity values only.
The functions $J_{H,E}=J_{H,E}(x,y,z,t)$ represent given source terms.
If, for the moment, we assume that the homogeneous Dirichlet boundary conditions 
are supplied with~\eqref{mxw}, then a standard Yee finite difference discretization
on a staggered Cartesian mesh results in a time-continuous space-discretized
system
\begin{equation}
\label{mxw0}
\begin{bmatrix}
M_\mu & 0 \\
0    & M_\eps    
\end{bmatrix}
\begin{bmatrix}
\bm{h}' \\ \bm{e}'  
\end{bmatrix} =
-\begin{bmatrix}
M_{\sigma_1} & K \\
-K^T       & M_{\sigma_2}    
\end{bmatrix}
\begin{bmatrix}
\bm{h} \\ \bm{e}  
\end{bmatrix}
+
\begin{bmatrix}
\bm{j}_H \\ \bm{j}_E  
\end{bmatrix},
\end{equation}
where the vector functions $\bm{h}(t)$ and $\bm{e}(t)$
contain the mesh values of the unknown fields,
$M_\mu$, $M_\eps$, and $M_{\sigma_{1,2}}$ are 
diagonal matrices containing
the values of $\mu$, $\eps$, and $\sigma_{1,2}$, respectively,
$K$ and $K^T$ are discrete curl operators and 
$\bm{j}_{H,E}(t)$ are the mesh values of the source
functions $J_{H,E}$.

Note that a semidiscrete system of ordinary differential
equations (ODEs), which is very similar to~\eqref{mxw0}, is
also obtained when the standard Whitney-N\'ed\'elec vector finite 
elements are employed (see 
e.g.~\cite{RodrigueWhite01,BotchevVerwer09,VerwerBotchev09}).
In this case $M_\mu$, $M_\eps$, and $M_{\sigma_{1,2}}$ are the mass
matrices.  It is convenient to rewrite the system~\eqref{mxw0}
as
\begin{equation}
\label{mxw1}
\begin{gathered}
\begin{bmatrix}
\bm{h}' \\ \bm{e}'  
\end{bmatrix} =
-A \cdot
\begin{bmatrix}
\bm{h} \\ \bm{e}  
\end{bmatrix}
+
\begin{bmatrix}
M_\mu^{-1}\bm{j}_H(t) \\ M_\eps^{-1}\bm{j}_E(t)  
\end{bmatrix},
\\
A = \begin{bmatrix}
M_\mu^{-1} & 0 \\
0    & M_\eps^{-1}    
\end{bmatrix}
\begin{bmatrix}
M_{\sigma_1} & K \\
-K^T       & M_{\sigma_2}    
\end{bmatrix}
=
\begin{bmatrix}
M_{1} & K_1 \\
-K_2^T       & M_{2}    
\end{bmatrix}
\in\Rr^{n\times n},
\end{gathered}
\end{equation}
where
$M_1=M_\mu^{-1}M_{\sigma_1}$,
$K_1=M_\mu^{-1}K$,
$M_2=M_\eps^{-1}M_{\sigma_2}$,
$K_2^T=M_\eps^{-1}K^T$,
and the inverse mass matrices are computed
explicitly only if they are diagonal or block diagonal.
The latter is the case if discontinuous Galerkin
finite elements are used, see e.g.~\cite{Sarmany_ea2013}.
We denote the size of the ODE system in~\eqref{mxw1} by $n$,
and let $n=n_1+n_2$, where $n_{1,2}$ are the numbers of degrees
of freedom associated with magnetic and electric fields,
respectively.
 
Employment of the nonreflective PML boundary conditions~\cite{Johnson2010_PML,PML94} 
means that auxiliary variables are added to the Maxwell system~\eqref{mxw}
which, after space discretization, enter the semidiscrete system~\eqref{mxw1}
as well.
These additional variables are nonzero only in the so-called PML region 
(a region just outside the boundary of the domain of interest). 
Incorporation of the PML boundary conditions into~\eqref{mxw},\eqref{mxw1}
(for a detailed derivation we refer to~\cite{deCloetMarissenWestendorp2015})
leads to the resulting semi-discrete ODE system of an extended size
\begin{equation}
\label{mxw2}
y'(t)=-\Aa y(t) + g(t),
\qquad
\Aa = 
\begin{bmatrix}
A & B_1^T \\ -B_2 & 0  
\end{bmatrix}\in\Rr^{N\times N},
\end{equation}
where $N=m+n$, with $m$ being the number of space-discretized auxiliary PML 
variables ($m$ is proportional to the number of mesh points in the PML region),
and the matrices $B_{1,2}\in\Rr^{m\times n}$ couple these variables to the 
main variables $\bm{h}$ and $\bm{e}$.
For representative values
of $m$ and $n$ see Table~\ref{t:mesh} in the numerical experiment
section below.
The matrices $B_{1,2}$ are defined in more detail below in~\ref{app1}.

\section{Nested Schur complement solver}
Exponential time integration based on rational 
shift-and-invert Krylov subspace
methods~\cite{EshofHochbruck06,Botchev2016} as well as 
implicit time integration~\cite{VerwerBotchev09} of 
systems~\eqref{mxw2} involves solution of linear systems
\begin{equation}
\label{ls}
(I+\gamma\Aa) x = b,
\end{equation}
where $\Aa$ is defined in~\eqref{mxw2} and $\gamma>0$ is a given parameter, 
related (or equal) to the time step size.
Matrices having nested saddle point structure\footnote{Formally speaking,
  $I+\gamma\Aa$ gets a saddle point structure when we switch the sign in the second
  block row of a linear system $(I+\gamma\Aa) x =b$.} as the matrix $\Aa$ are recently
called double saddle point problems~\cite{BeikBenzi2018}.
Our starting point in construction of preconditioners for matrices of this type
is the observation that an efficient preconditioner should involve 
a Schur complement 
(see e.g.~\cite{MurphyGolubWathen2000,ActaNumerKKT2005,ElmanSilvesterWathen:book}).
In particular, for matrices
$$
\begin{bmatrix}
\At & \Bt \\ \Ct & \Dt  
\end{bmatrix},
$$
good Schur complement-based block-diagonal preconditioners are
$$
\begin{bmatrix}
\At   & \Ot \\
\Ot   & \Dt-\Ct\At^{-1}\Bt   
\end{bmatrix}
\quad\text{or}\quad
\begin{bmatrix}
\At-\Bt\Dt^{-1}\Ct & \Ot \\
\Ot               & \Dt
\end{bmatrix}.
$$
For modern Krylov subspace methods such as GMRES, these preconditioners guarantee
convergence in three iterations~\cite{MurphyGolubWathen2000}.
Applying them for our matrix $I+\gamma\Aa$ means that linear systems with
$$
\begin{aligned}
\text{either (option~1)}\quad & \Dt-\Ct\At^{-1}\Bt= I+\gamma^2B_2(I+\gamma A)^{-1}B_1^T
\\
    \text{or (option~2)}\quad & \At-\Bt\Dt^{-1}\Ct= I+\gamma A+\gamma^2B_1^TB_2
\end{aligned}
$$
have to be solved efficiently.  Comparing these two possible options,
we choose option~2 because of the simpler structure of the matrix.
Furthermore, assuming for the moment that the systems with the
matrix $I+\gamma A$ can be solved efficiently and taking into
account that $B_1^TB_2$ is of a low rank, we may expect that 
$I+\gamma A$ can be a good preconditioner when solving systems
with the matrix $I+\gamma A+\gamma^2B_1^TB_2$.
This expectation is confirmed in practice: 
number of preconditioned by $I+\gamma A$ iterations to solve 
systems with $I+\gamma A+\gamma^2B_1^TB_2$ remain approximately 
constant as the spatial discretization mesh gets finer 
(see Table~\ref{t:B1B2}).
Moreover, as shown by formula~\eqref{B1B2} in the appendix below,
the matrix $\gamma^2B_1^TB_2$ depends on the mesh size in
a way similar to the matrix $I+\gamma A$ does: only
the $(1,2)$ and $(2,1)$ blocks in this matrix 
depend on the mesh size
as $\sim 1/h$ (assuming $h=h_x=h_y=h_z$). 

Now a question arises whether and how the systems with the matrix
\begin{equation}
\label{I_gA}
I+\gamma A =
\begin{bmatrix}
I + \gamma M_1 & \gamma K_1 \\ -\gamma K_2^T & I+\gamma M_2  
\end{bmatrix}
\end{equation}
can be solved efficiently.  We proceed in a similar way
as for the matrix $I+\gamma\Aa$ and explore the two possible options 
for a Schur complement-based preconditioner:
\begin{align}
\label{sch1}
\text{option 1:}\quad & I+\gamma M_2 + \gamma^2K_2^T(I+\gamma M_1)^{-1}K_1,
\\
\text{option 2:}\quad & I+\gamma M_1 + \gamma^2K_1(I+\gamma M_2)^{-1}K_2^T.
\notag
\end{align}
In many applications involving the Maxwell equations (such as, e.g.,
photonics and gas-and-oil exploration) the permeability
$\mu$ is usually constant ($\mu\equiv 1$), 
whereas $\eps$ is not and can be a strongly varying function.
Since the motivation for this work is photonics modeling,
we choose for option~1, where the matrix
$\gamma^2K_2^T(I+\gamma M_1)^{-1}K_1$ has a simpler structure than
the matrix $\gamma^2K_1(I+\gamma M_2)^{-1}K_2^T$ in option~2.
It is convenient to rewrite the chosen Schur complement~\eqref{sch1}
in the form
\begin{equation}
\label{sch2}
I+\gamma M_2 + \gamma^2K_2^T(I+\gamma M_1)^{-1}K_1 =
M_\eps^{-1}
\left[ M_\eps +\gamma M_{\sigma_2} + \gamma^2 K^T(M_\mu + \gamma M_{\sigma_1})^{-1}K
\right],
\end{equation}
which has an advantage that the matrix in brackets is symmetric
positive definite.
Further inspection of the bracketed matrix reveals that 
it is similar in structure to a shifted Laplacian, where 
the shift is given by $M_\eps +\gamma M_{\sigma_2}$.
For this reason a large variety of solvers is available for
solving systems with this matrix.  These include
(i)~sparse direct factorization solvers (on coarse meshes);
(ii)~multigrid solvers (which should not be too difficult in
implementation since only one field is involved);
(iii)~algebraic multigrid methods;
(iv)~preconditioned conjugate gradients (CG).
In this paper we use the CG solver preconditioned by
the incomplete Cholesky IC(0) preconditioner~\cite{ICCG}.

The described nested Schur complement approach can be used
either as a Schur complement based preconditioner or 
as a ``direct'' solver, computing the inverse 
action $(I+\gamma\Aa)^{-1}$ with the inner two-level iterative
solver for the Schur complement.
In Figure~\ref{f:alg} we outline the outer part of the 
introduced nested Schur complement solver for the case 
the action of $(I+\gamma\Aa)^{-1}$ is computed.
The inner part, corresponding to the action of 
$(I+\gamma A)^{-1}$, can then be computed in a similar
fashion, as described above.

\begin{figure}
\centering{\begin{minipage}{0.8\linewidth}
\fbox{Given $I+\gamma\Aa\in\Rr^{N\times N}$ and $b\in\Rr^N$
(cf.~\eqref{mxw1},\eqref{mxw2}),
solve $(I+\gamma\Aa)x=b$}
\\
1.~Partition $b$ into 
$b=\begin{bmatrix}b_1\\b_2\end{bmatrix}$, 
with $b_1\in\Rr^n$ and $b_2\in\Rr^m$.
\\
2.~Set $b_1:=b_1-\gamma B_1^Tb_2$.
\\
3.~The outer solver: solve $(I+\gamma A + \gamma^2 B_1^TB_2)x_1=b_1$
iteratively,\\$\phantom{\text{3.~}}$preconditioned by $I+\gamma A$.
\\
4.~Set $x_2:=b_2 + \gamma B_2x_1$  and 
$x:=\begin{bmatrix}x_1\\x_2\end{bmatrix}$.
\end{minipage}}
\caption{An algorithm description for the outer part of the nested 
Schur complement solver}
\label{f:alg}  
\end{figure}

\section{Other possible preconditioners}
A variety of other different preconditioners are available for 
solving saddle point problems of type~\eqref{ls}, see 
e.g.~\cite{ActaNumerKKT2005,ElmanSilvesterWathen:book}
and recent work~\cite{BeikBenzi2018}.
However, a general problem with linear solvers employed
in implicit and exponential time integrators is that 
the additional computational work spent for solving linear systems
has to be paid off by an increase in a time step.
Assume, for instance, that approximately ten iterations with a 
basic preconditioner have to be done per time step,
such that costs of a preconditioned matrix--vector product (matvec)
are approximately equal to the costs of an unpreconditioned
matvec (which can be achieved by the Eisenstat's 
trick~\cite{EisenTrick}).
Then the time step size 
have to be increased at least by a factor of ten to compensate
for the increased costs.
Such a time step increase, however, is not always possible due to 
accuracy restrictions, especially
for mildly stiff problems such as the Maxwell equations with PML
boundary conditions.
This makes the choice of a proper preconditioner difficult and
significantly restricts a variety of possible 
options~\cite{VerwerBotchev09}.

\begin{figure}
\begin{center}
\includegraphics[width=0.8\linewidth]{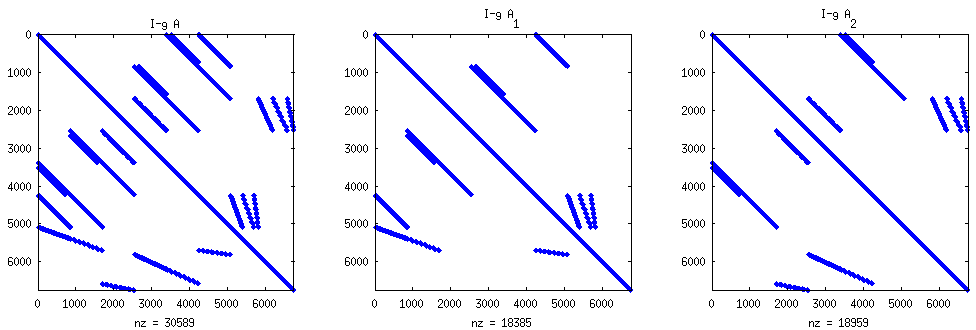}
\\
\includegraphics[width=0.8\linewidth]{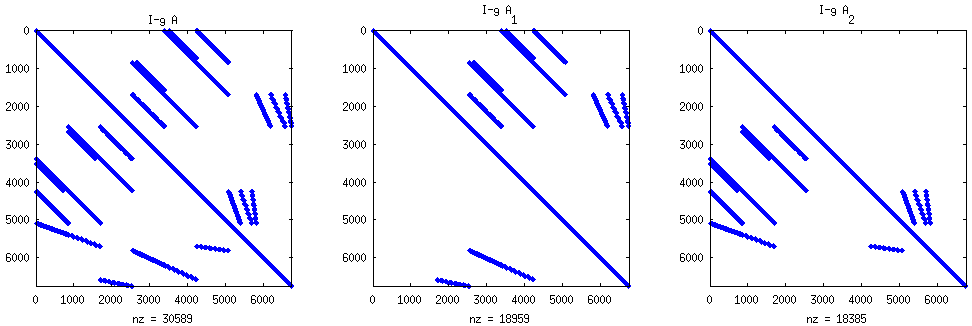}
\end{center}
\caption{Sparsity patterns of the matrix $I+\gamma A$ (left column)
  and its factors $I+\gamma A_{1,2}$ in the
  ADI (top row) and FS (bottom row)
  preconditioners for a coarse mesh $10\times 10\times 6$.
  The matrices $I+\gamma A$ and $I+\gamma A_{1,2}$ occupy the first 5082 rows
  and columns of the matrices $I+\gamma\Aa$ and $I+\gamma\Aa_{1,2}$,
  respectively.}
\label{f:prec}   
\end{figure}
Recently an efficient alternating direction implicit (ADI)
preconditioner is proposed and analyzed for solving time-dependent Maxwell 
equations~\cite{Hochbruck_ea2015_ADI} discretized in space by finite
differences. 
In~\cite{deCloetMarissenWestendorp2015} de~Cloet, Marissen and Westendorp
compared performance of this ADI preconditioner with another preconditioner
based on the field splitting (i.e., a splitting into the magnetic
and electric fields).  Their conclusion is that this field splitting (FS)
preconditioner outperforms the ADI preconditioner in terms of the CPU time.
Unlike the ADI preconditioner, the FS preconditioner is not restricted 
to the finite difference approximations on Cartesian meshes.
The linear system~\eqref{ls}, with either the FS or ADI preconditioner
applied from the right, can be written as
\begin{equation}
\label{pls}
\widetilde{\Aa}\widetilde{x} = b,
\quad \widetilde{\Aa}=(I+\gamma\Aa)\Mm^{-1},
\quad \widetilde{x} = \Mm x,
\end{equation}
where
\begin{equation}
\label{fsp}
\Mm=(I+\gamma\Aa_1)(I+\gamma\Aa_2), \quad \Aa_1+\Aa_2 = \Aa.
\end{equation}
In the FS preconditioner the matrices $\Aa_{1,2}$ are defined as 
\begin{equation}
\label{FS}
\begin{gathered}
\Aa_1 =
\begin{bmatrix}
A_1 & B_{1,H}^T \\ -B_{2,H} & 0  
\end{bmatrix},
\quad
\Aa_2 =
\begin{bmatrix}
A_2 & B_{1,E}^T \\ -B_{2,E} & 0  
\end{bmatrix},
\\
A = A_1 + A_2, \quad
A_1 =
\begin{bmatrix}
M_1 & K_1 \\ 0 & 0  
\end{bmatrix},
\quad
A_2 =
\begin{bmatrix}
0   &  0 \\ -K_2^T   & M_2
\end{bmatrix},
\end{gathered}
\end{equation}
where the matrices $B_{j,H}$, $B_{j,E}$, $j=1,2$, form splittings
of the PML blocks $B_{1,2}$,
$$
B_1=B_{1,H}+ B_{1,E}, \quad
B_2=B_{2,H}+ B_{2,E},
$$
defined in~\ref{app2}.
For complete definition of the ADI preconditioner we refer to
\cite[Section~5.1]{deCloetMarissenWestendorp2015} and
\cite{Hochbruck_ea2015_ADI}.
In both the ADI and FS preconditioners the factors
$I+\gamma\Aa_{1,2}$ are not triangular matrices (see Figure~\ref{f:prec})
and sparse LU~factorizations have to be carried out to implement
the preconditioner actions.

\begin{figure}
\centering{%
\includegraphics[width=0.35\linewidth]{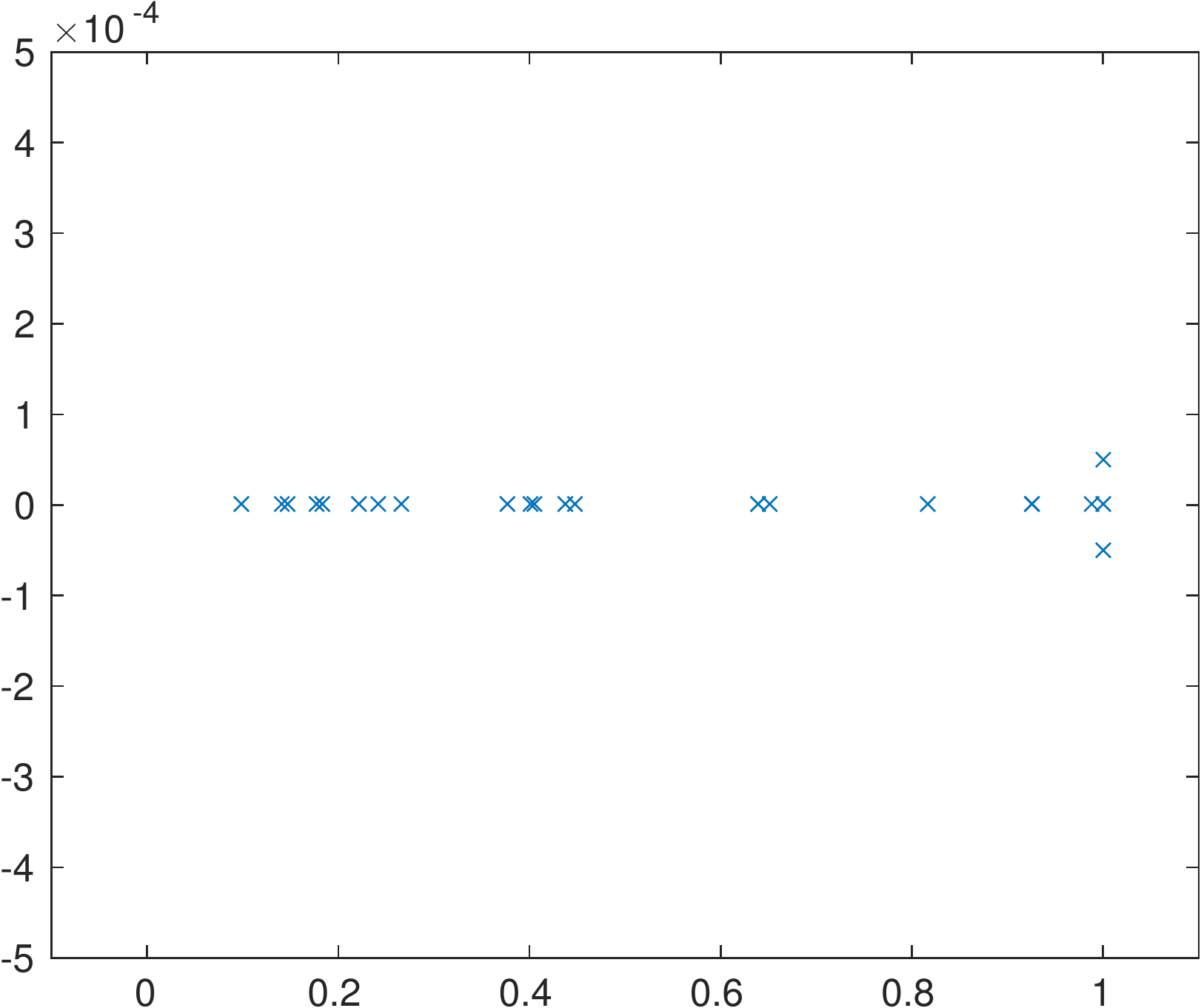} ~~
\includegraphics[width=0.35\linewidth]{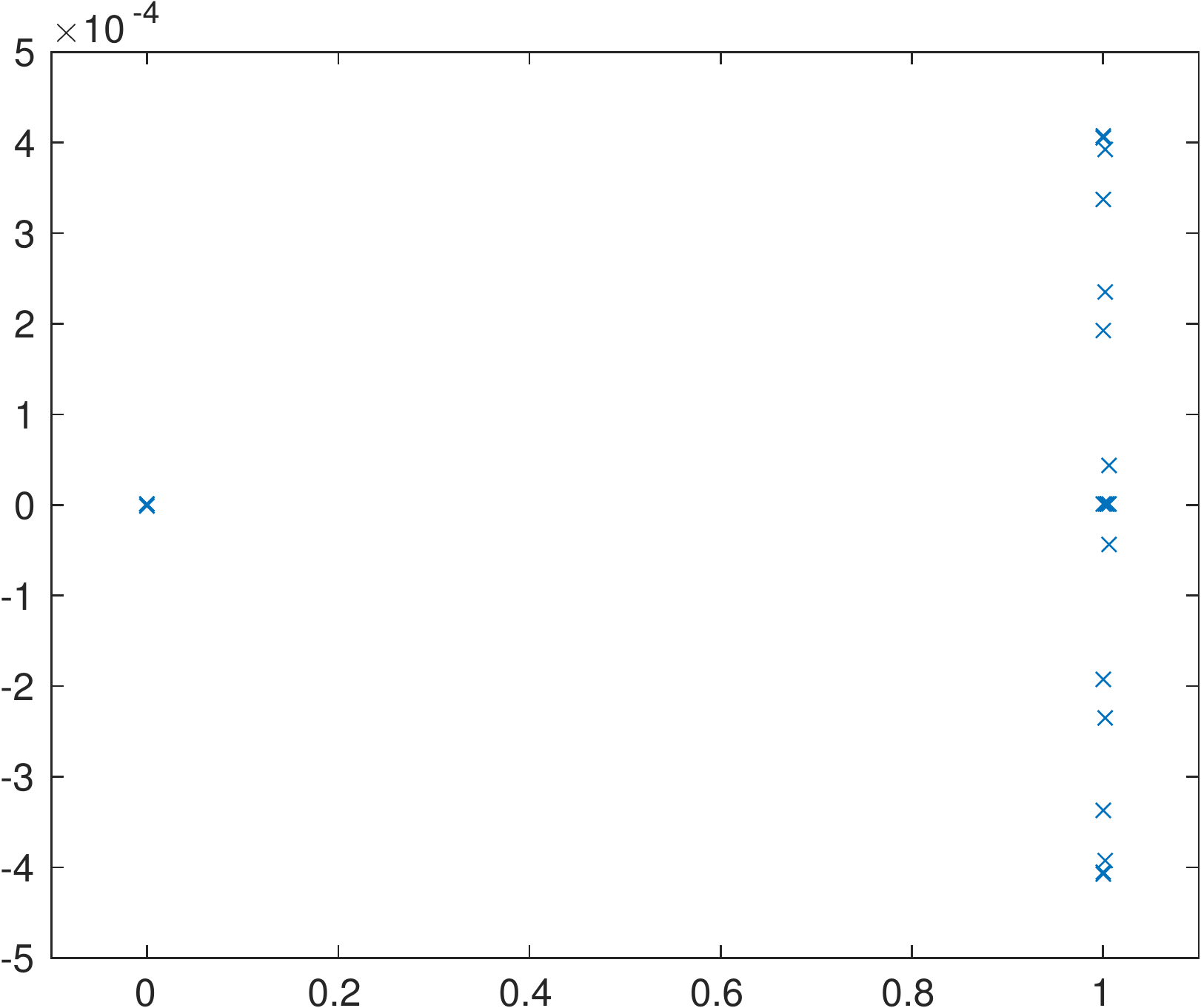}}
\caption{Ritz values on the complex plane computed after 
24~iterations of the FOM iterative method
with the ADI (left) and FS (right) preconditioners}
\label{Rvv}
\end{figure}

Numerical tests presented in~\cite{deCloetMarissenWestendorp2015}
demonstrate that both ADI and FS preconditioners require approximately
the same number of iterations to converge.  However, the FS preconditioner
is faster than ADI in terms of the CPU time.  This is caused by another
attractive property of the FS preconditioner observed
in~\cite{deCloetMarissenWestendorp2015}: sparse LU factorizations
of the FS factors $I+\gamma\Aa_{1,2}$ do not yield any additional fill in.
This is the case for the ADI preconditioner,
see \cite[Section~6.2]{deCloetMarissenWestendorp2015},
where the fill in varies from 25\% on coarse meshes
to 65\% on the mesh $80\times 80\times 48$.
Ritz values obtained after 24~iterations of the FOM iterative
method~\cite{GMRES} with both the ADI and FS preconditioners
are plot in Figure~\ref{Rvv}.  There we see that both preconditioners
yield preconditioned matrices with effectively real spectra.

\section{Numerical experiments}
In the test problem a 3D photonic crystal
is considered.
At the $x$- and $y$- boundaries of the spatial domain $[1,4]\times [1,4]\times [0,3]$
the PML boundary conditions are imposed,
whereas homogeneous Dirichlet (perfectly conducting) boundary
conditions are posed on the $z$-boundaries.
The PML regions extend
the total computational domain along the $x$- and $y$-walls to 
$[0,5]\times [0,5]\times [0,3]$.
The crystal consists of $3\times 3\times 3$
spheres of radius $0.4$ centered at points
$(x_i,y_j,z_k)=(2.5+i,2.5+j,1.5+k)$, $i,j,k=-1,0,1$.
The magnetic permeability $\mu\equiv 1$ in the whole domain,
whereas the electric permittivity $\eps$ is set to~$8.9$ inside 
the spheres and to~$1$ everywhere else in the domain.

\begin{table}
\caption{The number of degrees of freedom for the meshes used in the tests}
\label{t:mesh}
\centering\begin{tabular}{ccc}
\hline\hline
mesh                       & system size 
\\ 
$n_x \times n_y \times n_z$ & $N=n+m$     
\\
\hline
$20\times 20\times 12$      & $ 45\,565= 34\,398+11\,167$           
\\
$40\times 40\times 24$      & $333\,425=252\,150+81\,275$           
\\
$80\times 80\times 48$      & $2\,548\,441=1\,928\,934+619\,507$    
\\
$160\times 160\times 96$    & $ 19\,922\,345 = 15\,086\,022 + 4\,836\,323$
\\\hline
\end{tabular}
\end{table}

We consider matrices $\Aa$ resulting from the standard Yee finite difference
approximation, see Table~\ref{t:mesh} for typical
mesh sizes used in the tests.
The parameter $\gamma$ is chosen as explained in~\cite{Botchev2016}
and set to $\gamma=0.012$ in all the tests.
The tests are run in Matlab on a Linux PC
with 8~CPUs Intel Xeon~E5504 2.00GHz.  

\begin{table}
\caption{Iteration numbers and residual norms for 
solving linear systems with $I+\gamma A + \gamma^2 B_1^TB_2$ 
preconditioned by $I+\gamma A$ and norms of 
the symmetric and skew-symmetric parts
of $\gamma^2B_1^TB_2$, 
$\Ht=\frac{\gamma^2}2(B_1^TB_2+(B_1^TB_2)^T)$
$\St=\frac{\gamma^2}2(B_1^TB_2-(B_1^TB_2)^T)$.}
\label{t:B1B2}
\centerline{\begin{tabular}{cccc}
\hline\hline
mesh                    & \# iter,  & $\|\Ht\|_1 $      & $\|\St\|_1$     \\
                        & resid.norm&                   &           
\\\hline                         
$20\times 20 \times 12$ & 21, \texttt{2.80e-07} & 1968  & 9.8              
\\                               
$40\times 40 \times 24$ & 21, \texttt{4.73e-07} & 1572  & 19.6              
\\                               
$60\times 60 \times 36$ & 21, \texttt{5.71e-07} & 1458 & 29.4              
\\
\hline
\end{tabular}}
\end{table}

\begin{figure}
\centering{\includegraphics[width=0.7\linewidth]{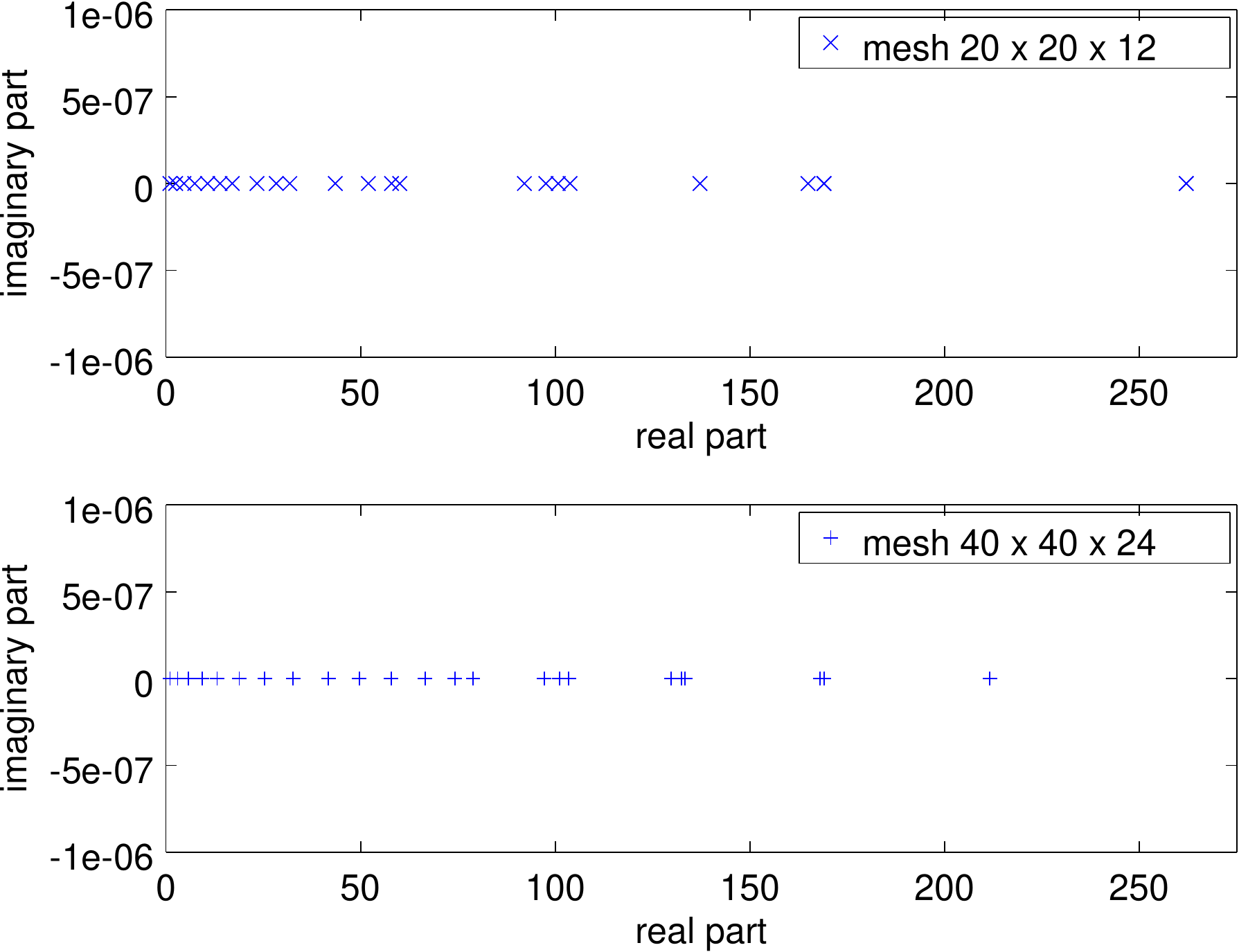}}
\caption{Ritz values of the preconditioned matrix $I+(I+\gamma A)^{-1}\gamma^2 B_1^TB_2$
on the complex plane for 
the mesh $20\times 20\times 12$ (top) and 
$40\times 40\times 24$ (bottom)}
\label{f:I_gA}  
\end{figure}

In Table~\ref{t:B1B2} we illustrate the fact that linear systems
with the matrix $I+\gamma A + \gamma^2 B_1^TB_2$ can be efficiently
solved iteratively using $I+\gamma A$ as a preconditioner.
The preconditioner actions here are carried with the help of
the sparse LU factorization (UMFPACK in Matlab), that is why in 
this case we can not use a fine mesh.
In this test the exact solution vectors are taken to have normally distributed 
random entries with zero mean and variance one.  This is done
to make the test difficult so that the solver can not possibly profit 
from the solution smoothness.
The number of iterations listed there are for the BiCGstab(2)
iterative solver (the standard built-in solver in Matlab)
run to satisfy the stopping criterion tolerance of $10^{-6}$. 
In the same table we also list the norms of the symmetric and 
skew-symmetric parts of $\gamma^2 B_1^TB_2$.  The values 
show that the field of values of this matrix is bounded and 
confirm
the expectation given by relation~\eqref{B1B2}: only the off-diagonal
blocks (related to the skew-symmetric part) of this matrix depend 
on the mesh size and this dependence is linear.
In Figure~\ref{f:I_gA} we plot 24~Ritz values of the preconditioned
matrix $I+(I+\gamma A)^{-1}\gamma^2 B_1^TB_2$.  As we see, the Ritz
values are real and well clustered, which means that the preconditioner
is efficient and damps the skew-symmetric part of the system matrix
well.

\begin{table}
\caption{Results for the nested Schur complement solver for the ``difficult'' 
test case, with random exact solution vector.}
\label{t:NSC}
\centering{\begin{tabular}{cccc}
\hline\hline
mesh                   & residual           & iterations   & CPU    \\
                       & norm               & outer (inner)& time   \\\hline
$20\times 20\times 12$ & \texttt{7.50e-11}  & 31 (68)    & 3.59~s 
\\
$40\times 40\times 24$ & \texttt{6.67e-11}  & 32 (108)   & 17.3~s 
\\
$80\times 80\times 48$ & \texttt{6.84e-11}  & 32 (145)   & 135.6~s 
\\
$160\times 160\times 96$
                       & \texttt{7.01e-11}   & 31 (234)   & 1258~s 
\\\hline
\end{tabular}}
\end{table}

We now test our nested Schur complement solver. 
The linear systems in the test have the exact solutions
which are again a normally distributed 
random vector whose entries have zero mean and variance one.
This is done to prevent the solver from profiting
from the solution smoothness.
The solver is applied in its direct form, as outlined in Figure~\ref{f:alg},
with the outer solver GMRES(10) and inner solver ICCG(0).
The stopping criterion tolerance in both solvers is set to $10^{-10}$.
We see that the number of the outer iterations remains constant, independently 
of the mesh size, as expected.
The number of inner iterations grows because the CG solver is used
with the simple IC(0) preconditioner.
We note that the number of inner iterations changes from one outer iteration
to another and the inner iteration count reported in the table is the maximum
number of inner iterations (required in all cases at the last outer iteration).

\begin{table}
\caption{Comparison results of the nested Schur complement solver 
and the FS preconditioner.  The former is employed with GMRES(10) and 
ICCG(0) as the outer and inner solvers, respectively.  The FS preconditioner
is applied with nonrestarted GMRES.}
\label{t:vs}
\centering{\begin{tabular}{cclc}
  \hline\hline
  mesh,                  & method & CPU time,    &  residual   \\
  tolerance              &        & iter outer(inner) &  norm       \\
  \hline
  $40\times 40\times 24$ & FS prec.  & 0.49~s, 7~~~& \texttt{7.58e-12}\\ 
  \texttt{9.64e-11}      & nested Schur & 0.54~s, 2(8)& \texttt{1.23e-13}\\
  \hline
  $80\times 80\times 48$ & FS prec.  & 3.84~s, 8~~~& \texttt{5.95e-10}\\ 
  \texttt{8.09e-09}      & nested Schur & 3.51~s, 2(8)& \texttt{2.09e-09}\\
  \hline                                               
  $160\times160\times 96$& FS prec.  & 59.7~s, 14~~~&\texttt{3.97e-09} \\
  \texttt{4.23e-09}      & nested Schur & 45.8~s, 2(19)&\texttt{3.77e-09}\\
  \hline
\end{tabular}}
\end{table}

Finally, in Table~\ref{t:vs} we present results of comparison of the nested
Schur complement solver with the FS preconditioner.  The comparison is done
on linear systems arising in the time integration carried out by 
an exponential integrator based on a rational shift-and-invert exponential 
Krylov subspace method.  Therefore the stopping criterion tolerance in the tests 
vary and is reported in the table.
We see that the nested Schur complement solver outperforms the FS preconditioner
on fine meshes.  This is expected because the FS preconditioner does not
converge mesh-independently.
 
\section{Conclusions and an outlook to further research}
A nested Schur complement solver is proposed for iterative linear 
system solution within exponential and implicit time integration 
of the Maxwell equations.  The solver exhibits a mesh-independent 
convergence and outperforms other preconditioners, such 
ADI (alternative direction implicit) and FS (field splitting)
preconditioners.

Different aspects in the proposed concept require further
investigation and possible improvement.
In future we plan to test the nested Schur complement solver
in combination with a more robust (and mesh-independent) 
inner iterative solver.  
Another interesting research
question is which form of the solver is most efficient:
its direct form, as tested in this paper, or iterative,
as a three-level iterative solver.

\bibliography{matfun,mxw,my_bib}
\bibliographystyle{abbrv}

\appendix

\section{PML matrices $B_{1,2}$}
\label{app1}
For simplicity of the presentation we assume that 
the Yee finite difference discretization is used on
a mesh of $n_x\times n_y\times n_z$ cells.  Then it is 
possible to construct this discretization in such a way that 
all the three components of each of the two fields have
the same number of degrees of freedom $(n_x+1)(n_y+1)(n_z+1)$,
This can be done by augmenting the ``shorter'' components 
with auxiliary ``void'' degrees of freedom.
Then $n_1=n_2=n/2$ and the matrices $B_{1,2}$ are defined with the help of 
the following matrices $\widehat{B}_{1,2}\in\Rr^{2n\times n}$:
\begin{equation}
\label{B12}
\widehat{B}_1^T=\begin{bmatrix}
K_1 & 0      & -I &  0 \\
0   & -K_2^T & 0  & -I 
\end{bmatrix},
\quad
\widehat{B}_2 = 
\begin{bmatrix}
0      & \Sigma \\
\Sigma &  0 \\
-\Sigma_*       & 0 \\
0         & -\Sigma_*  
\end{bmatrix},
\end{equation}
where diagonal matrices $\Sigma,\Sigma_*\in\Rr^{n_1\times n_1}$ read
$$
\Sigma = \mathrm{diag} 
\begin{bmatrix}
\bm{\sigma}_x ,& \bm{\sigma}_y ,& \bm{\sigma}_z
\end{bmatrix},
\qquad
\Sigma_* = \mathrm{diag} 
\begin{bmatrix}
\bm{\sigma}_y\cdot\bm{\sigma}_z ,& 
\bm{\sigma}_x\cdot\bm{\sigma}_z ,& 
\bm{\sigma}_x\cdot\bm{\sigma}_y
\end{bmatrix}.
$$
Here $\bm{\sigma}_{x,y,z}$ are the PML artificial conductivity 
values for the corresponding direction and $\cdot$ denotes
elementwise multiplication.  Since $\bm{\sigma}_{x,y,z}$
are only nonzero inside the PML region, there are a lot of
redundant zero rows in $\widehat{B}_{1,2}$.  By omitting these
zero rows we obtain the matrices $B_{1,2}$.
Note also that
\begin{equation}
\label{B1B2}
\widehat{B}_1^T\widehat{B}_2 = B_1^TB_2 =
\begin{bmatrix}
\Sigma_*     & K_1\Sigma \\
-K_2^T\Sigma & \Sigma_* 
\end{bmatrix},
\end{equation}
which means that only the $(1,2)$ and $(2,1)$ blocks in this matrix 
depend on the mesh (as $K_{1,2}$ are discrete curl operators)
as $\sim 1/h$ (assuming $h=h_x=h_y=h_z$).
This is important in construction of preconditioners
for linear system solution with the matrix $I+\gamma A + \gamma^2 B_1^TB_2$.

\section{Field splitting of the matrices $B_{1,2}$}
\label{app2}
To define the FS preconditioner we split the matrices
$\widehat{B}_{1,2}$ into the components corresponding to
the two fields $H$ and $E$ as follows:
$$
\widehat{B}_1 = \widehat{B}_{1,H} + \widehat{B}_{1,E}
\qquad
\widehat{B}_2 = \widehat{B}_{2,H} + \widehat{B}_{2,E},
$$
where
\begin{equation}
\label{B12H}
\widehat{B}_{1,H}^T:=\begin{bmatrix}
K_1 & 0      & -I &  0 \\
0   & 0 & 0  & 0
\end{bmatrix},
\quad
\widehat{B}_{2,H} := 
\begin{bmatrix}
0          & \Sigma \\
0          &  0 \\
-\Sigma_*  & 0 \\
0          & 0   
\end{bmatrix},
\end{equation}
and
$\widehat{B}_{1,E}:=\widehat{B}_1-\widehat{B}_{1,H}$,
$\widehat{B}_{2,E}:=\widehat{B}_2-\widehat{B}_{2,H}$.
The matrices
${B}_{1,H}$,
${B}_{2,H}$,
${B}_{1,E}$,
${B}_{2,E}$
are then obtained from 
the matrices 
$\widehat{B}_{1,H}$,
$\widehat{B}_{2,H}$,
$\widehat{B}_{1,E}$,
$\widehat{B}_{2,E}$,
respectively, 
by omitting zero rows in the latter four.

\end{document}